\documentclass{amsart}
\usepackage[T1]{fontenc}
\usepackage{ae}
\usepackage{aecompl}%
\usepackage{ucs}
\usepackage[utf8x]{inputenc}
\usepackage{color}
\usepackage[enableskew,vcentermath]{youngtab}
\usepackage{amsfonts}
\usepackage{hhline}
\usepackage{amsthm}
\usepackage{amsmath}
\usepackage{amssymb}
\usepackage[all]{pstricks}
\usepackage[justification=centering,labelfont=bf]{caption}
\usepackage{hyperref}

\begin{document}

\numberwithin{equation}{section}

\newtheorem{Le}{Lemma}[section]
\newtheorem{Ko}[Le]{Lemma}
\newtheorem{Sa}[Le]{Theorem}
\newtheorem{pro}[Le]{Proposition}

\newtheorem{theorem}[Le]{Theorem}
\newtheorem{lemma}[Le]{Lemma}
\newtheorem{Con}[Le]{Conjecture}
\newtheoremstyle{Bemerkung}
  {}{}{}{}{\bfseries}{.}{0.5em}{{\thmname{#1}\thmnumber{ #2}\thmnote{ (#3)}}}

\theoremstyle{Bemerkung}
\newtheorem{definition}[Le]{Definition}
\newtheorem{Def}[Le]{Definition}
\newtheorem{example}[Le]{Example}

\newtheorem{remark}[Le]{Remark}
\newtheorem{Bem}[Le]{Remark}
\newtheorem{Bsp}[Le]{Example}

\renewcommand{\l}{\lambda}
\renewcommand{\L}{\lambda}
\newcommand{\bl}{{\bar\lambda}}
\newcommand{\bn}{\bar\nu}
\newcommand{\bp}{\bar p}
\newcommand{\A}{\mathcal{A}}
\newcommand{\B}{\mathcal{B}}
\newcommand{\C}{\mathcal{C}}
\newcommand{\D}{\mathcal{D}}
\renewcommand{\S}{\mathcal{S}}
\renewcommand{\L}{\mathcal{L}}
\newcommand{\E}{\mathcal{E}}
\newcommand{\N}{\mathcal{N}}
\renewcommand{\a}{\alpha}
\renewcommand{\b}{\beta}
\renewcommand{\c}{\gamma}
\renewcommand{\d}{\delta}
\newcommand{\e}{\epsilon}
\newcommand{\h}{\hfil}
\newcommand{\X}{X}
\newcommand{\abs}[1]{\left| #1 \right|}
\newcommand{\lm}{\l/\mu}
\renewcommand{\ln}{\l/\nu}
\newcommand{\ab}{\a/\b}
\newcommand{\m}{\mu}
\newcommand{\n}{\nu}
\newcommand{\dn}{\d_n/\d_{n-1}}
\newcommand{\co}{<\cdot}
\renewcommand{\r}{\rho}

\renewcommand{\i}{i(\lm)}
\newcommand{\ih}{ih(\lm)}
\newcommand{\iv}{iv(\lm)}
\renewcommand{\o}{o(\lm)}
\newcommand{\oh}{oh(\lm)}
\newcommand{\ov}{ov(\lm)}
\renewcommand{\v}{v(\lm)}
\newcommand{\cc}{cc(\lm)}

\title{The skew diagram poset and components of skew characters}
\author[C. Gutschwager]{Christian Gutschwager}
\address{Institut für Algebra, Zahlentheorie und Diskrete Mathematik, Leibniz Universität Hannover,  Welfengarten 1, D-30167 Hannover}
\email{gutschwager (at) math (dot) uni-hannover (dot) de}


\begin{abstract}
 We investigate the poset of skew diagrams ordered by adding or forming the union of skew diagrams. We will show that a skew diagram which has at least $n$ convex corners to the upper left and also to the lower right is larger than the skew diagram consisting of $n$ disconnected single boxes. Using this property, we obtain lower bounds for the number of components, constituents and pairs of components which differ by one box in a given skew character.
\end{abstract}

\maketitle

\section{Introduction and Notation}
Characters of the symmetric group are being investigated since the beginning of the $20$th century. Skew characters of the symmetric group decompose in the same way as skew Schur functions and their decomposition corresponds to the decomposition of products of Schubert classes (see \cite{Gut}) as well as the decomposition of the restriction of irreducible affine Hecke algebras to the Iwahori-Hecke algebras (see \cite{Ram}). 

We introduce a new poset on the set of skew diagrams (Section~\ref{sec:poset}) which allows us to obtain results about skew characters (Section~\ref{sec:fcmain:bounds}).

We mostly follow the standard notation in \cite{Sag} or \cite{Stanley}. A partition $\l=(\l_1,\l_2,\ldots,\l_l)$ is a weakly decreasing sequence of non-negative integers where only finitely many of the $\l_i$ are positive. We regard two partitions as the same if they differ only by the number of trailing zeros and call the positive $\l_i$ the parts of $\l$. The length is the number of positive parts and we write $l(\l)=l$ for the length and $\abs{\l}=\sum_i \l_i$ for the sum of the parts. With a partition $\l$ we associate a diagram, which we also denote by $\l$, containing $\l_i$ left-justified boxes in the $i$-th row and we use matrix style coordinates to refer to the boxes.

 We  write $dp(\l)=n$ if the partition $\l$ has $n$ different parts. Furthermore we set
 \[\d_n=(n,n-1,n-2,\ldots,2,1).\]

The conjugate $\l^c$ of $\l$ is the diagram which has $\l_i$ boxes in the $i$-th column.

For $\m \subseteq \l$ we define the skew diagram $\lm$ as the difference of the diagrams $\l$ and $\m$ defined as the difference of the set of the boxes. Rotation of $\lm$ by $180^\circ$ yields a skew diagram $(\lm)^\circ$ which is well defined up to translation.

 A skew tableau $T$ is a skew diagram in which positive integers are written into the boxes. A semistandard tableau of shape $\lm$ is a filling of $\lm$ with positive integers such that the entries weakly increase amongst the rows from left to right and strictly increase amongst the columns from top to bottom. The content of a semistandard tableau $T$ is $\n=(\n_1,\ldots)$ if the number of occurrences of the entry $i$ in $T$ is $\n_i$. The reverse row word of a tableau $T$ is the sequence obtained by reading the entries of $T$ from right to left and top to bottom starting at the first row. Such a sequence is said to be a lattice word if for all $i,n \geq1$ the number of occurrences of $i$ among the first $n$ terms is at least the number of occurrences of $i+1$ among these terms. The Littlewood-Richardson (LR) coefficient $c(\l;\m,\n)$ equals the number of semistandard tableaux of shape $\lm$ with content $\n$ such that the reverse row word is a lattice word. We will call those tableaux LR tableaux. The LR coefficients play an important role in different contexts (see \cite{Sag} or \cite{Stanley} for further details).

A standard Young tableaux of shape $\l$ is a filling of $\l$ with the numbers $1,\ldots,\abs{\l}$ such that the entries increase in each row from left to right and in each column from top to bottom. The number of standard Young tableaux  of shape $\l$ is denoted by $f^\l$ which is given by the well known hook length formula 
\[f^\l=\frac{\abs{\l}!}{\prod( \textnormal{hook length})}.\] Obviously the number of standard Young tableaux with $n$ boxes $f_n$ is given by $f_n=\sum_{\l\vdash n} f^\l$. Notice furthermore that $f_n$ is also the number of involutions in the symmetric group $S_n$ plus $1$.

The irreducible characters $[\l]$ of the symmetric group $S_n$ are naturally labeled by partitions $\l\vdash n$. The skew character $[\lm]$ corresponding to a skew diagram $\lm$ is defined by the LR coefficients 
\[ [\lm]=\sum_\n c(\l;\m,\n) [\n]. \]

The translation symmetry gives $[\lm]=[\ab]$ if the skew diagrams of $\lm$ and $\ab$ are the same up to translation while rotation symmetry gives $[(\lm)^\circ]=[\lm]$.  The conjugation symmetry $c(\l^c;\m^c,\n^c)=c(\l;\m,\n)$ is also well known and furthermore we have $c(\l;\m,\n)=c(\l;\n,\m)$.

A basic skew diagram $\lm$ is a skew diagram which satisfies $\m_i<\l_i$ and $\m_i\leq \l_{i+1}$ for each $1\leq i \leq l(\l)$. This means that $\lm$ doesn't contain empty rows or column in $\lm$. Empty rows or columns of a skew diagram don't influence the filling and so deleting empty rows or columns doesn't change the skew character or LR fillings.

Let $\A$ and $\B$ be non-empty sub-diagrams of a skew diagram $\D$ such that the union of $\A$ and $\B$ is $\D$.
Then we say that the skew diagram $\D$ is \textit{disconnected} or \textit{decays} into the skew diagrams $\A$ and $\B$ if no box of $\A$ (viewed as boxes in $\D$) is in the same row or column as a box of $\B$. Notice, that this also covers the case, when $\B$ again decays into two subdiagrams and $\A$ is between those two. We write $\D=\A\otimes\B$ if up to translation $\D$ decays into $\A$ and $\B$ and normally write $\A$ and $\B$ as basic skew diagrams. A skew diagram is connected if it does not decay. If $\D=\A\otimes\B=\C$ then by translation symmetry $[\D]=[\C]$, so reordering $\A, \B$ doesn't change the skew character.

A skew character whose skew diagram $\D$ decays into the skew diagrams $\A,\B$ is equivalent to the product of the  characters of the disconnected diagrams induced to a larger symmetric group. We have
 \[[\D]=([\A]\times[\B])\uparrow_{S_n\times S_m}^{S_{n+m}}=:[\A]\otimes[\B]\]
 with $\abs{\A}=n,\abs{\B}=m$.  If $\D=\lm$ and $\A,\B$ are proper partitions $\a,\b$ then we have 
\[[\lm]= \sum_\n c(\l;\m,\n)[\n]=\sum_\n c(\n;\a,\b)[\n] =[\a]\otimes[\b].\]

\section{The poset of skew diagrams}\label{sec:poset}

For a partition $\l$ we can define a path starting to the right at the lower left corner of $\l$ and following the shape of $\l$ to the upper right corner, ending with an upward going segment. We write this path as a sequence of $v$'s and $h$'s denoting either a vertical or horizontal step.

For example the path corresponding to the partition $\l=(5,3,1,1,1)=\yng(5,3,1,1,1)$ is given by the sequence \[s=(hvvvhhvhhv).\]

For a skew diagram $\lm$ we define two paths, the outer and inner path. The outer path, whose sequence we denote by $o(\lm)$ or simply $o$, is the path of $\l$. The inner path, whose sequence we denote by $i(\lm)$ or $i$, starts in the lower left corner of $\lm$ upwards to the lower left corner of $\m$ (there is no upward step if $l(\l)=l(\m)$) follows the path of $\m$ and ends with steps to the right at the upper right, provided $\m_1<\l_1$.

So for $\lm=(5,3,1,1,1)/(5,2,1)$ (see Figure~\ref{pic0})
\begin{figure}[h]\begin{center}
\psset{xunit=0.35cm,yunit=0.35cm,runit=0.35cm}
\begin{pspicture*}(-.2,-.2)(10.2,10.2)
\psline(0,0)(2,0)(2,6)(6,6)(6,8)(10,8)(10,10)
\psline(0,0)(0,4)(2,4)
\psline(4,6)(4,8)(6,8)         
\psdot*(0,2)\psdot*(2,2)\psdot*(8,8)\psdot*(0,0)\psdot*(2,0)\psdot*(0,4)\psdot*(2,4)\psdot*(2,6)\psdot*(4,6)\psdot*(4,8)\psdot*(6,6)\psdot*(6,8) \psdot*(10,8)\psdot*(10,10)
\end{pspicture*}
\caption{$\lm=(5,3,1,1,1)/(5,2,1)$} \label{pic0}\end{center}
\end{figure}
\newline we have the paths $o=s$ as above and $i=(vvhvhvhhhv)$.

Let $\lm$ be a skew diagram. Then both the $o$ and $i$ sequence have $\l_1+l(\l)$ entries. Furthermore, both sequences have $\l_1$ entries $h$ and $l(\l)$ entries $v$. For all $j\leq \l_1+l(\l)$ the number of entries $h$ among the first $j$ entries of $i$ is at most the number of entries $h$  among the first $j$ entries of $o$. Otherwise the partition $\m$ wouldn't be contained in $\l$. Furthermore, if those numbers are equal then the $o$ and $i$ path touch each other after $j$ steps. From this follows, that for a basic skew diagram there is no $j< \l_1+l(\l)$ such that the number of entries $h$ among the first $j$ entries of $i$ and $o$ are the same and both sequences continue with the same entry in the $j+1$'st position.

The sum $\m+\n=\l$ of two partitions $\m,\n$ is defined by $\l_i=\m_i+\n_i$. The partition $\m\cup\n$ contains the parts of both $\m$ and $\n$. These operations are conjugate to another
\[(\m+\n)^c=\m^c\cup\n^c\] and do not commute
\[(\l\cup\m)+\n\neq(\l+\n)\cup\m. \]

For example, we have
\begin{align*}
 \yng(4,4,2,1)+\young(XXX,XX,X,X)&=\young(\h\X\h\h\h\X\X,\h\X\h\h\h\X,\h\X\h,\h\X), & \yng(4,3,2,2) \cup \young(XXXX,XX,X)&=\young(\h\h\h\h,XXXX,\h\h\h,\h\h,\h\h,XX,X).
\end{align*}

Note that $\m\cup\n=\m\cup\n_1\cup\n_2\cup\cdots\cup\n_{l(\n)}$. Because of this we sometimes say that for $\m+\n$ we insert the columns of $\n$ into $\m$ and for $\m\cup\n$ that we insert the rows of $\n$ into $\m$. Note that this $+$ and $\cup$ introduce a partial order on the set of partitions and we say that a partition $\l$ is larger than $\l'$ if $\l$ can be obtained from $\l'$ by repeatedly using the operations $+,\cup$ with arbitrary partitions in any order. This should not be confused with the lexicographic order. Note that the sequences of $\m$ and $\m+(1^n)$ differ by one $h$ if $n\leq l(\m)$. By symmetry the sequences of $\m$ and $\m\cup(n)$ differ by one $v$ if $n\leq \m_1$. To be more precise, the sequence of $\m+(1^n)$ ($n\leq l(\m)$) is obtained from $\m$ by inserting an  $h$ such that there are exactly $n$ entries $v$ to the right of the new $h$ which also means that there are exactly $l(\m)-n$ entries $v$ to the left.

For two skew diagrams $\A=\lm, \B=\l'/\m'$ we define the operations $\A+\B=\a/\b$ and $\A\cup\B=\a'/\b'$ by $\a=\l+\l',\b=\m+\m'$ and $\a'=\l\cup\l',\b'=\m\cup\m'$, respectively. Clearly $\A+\B$ and $\A\cup\B$ are then again skew diagrams. Usually we regard two skew diagrams as the same if they contain boxes in the same position but for this definition the underlying partitions $\l,\l',\m $ and $\m'$ are important because different choices for $\l$ and $\m$ would lead to different $\ab$. For example, we have $(2,1)/(1^2)=(2)/(1)=\young(\h)$. But if we add in both cases $(1^2)$ we would get $(3,2)/(1^2)=\young(\h\h,\h)\neq\young(:\h\h,\h)=(3,1)/(1)$. However, this will never cause any problem because in general we assume that the skew diagrams are basic.

On the set of basic skew diagrams we define a partial order as follows. Let $\A,\B$ be skew diagrams, then we say that $\A$ is greater or equal to $\B$ if there exists $n\in\mathbb{N}$ and for $1\leq i \leq n$ it is $\circ^i\in\{+,\cup\}$ and $\C^i$ a skew diagram so that we have
\[\A=\left(\cdots \left((\B\circ^1 \C^1)\circ^2\C^2\right)\cdots \right) \circ^n\C^n.\] 
Notice that it is not enough that $\l$ or $\m$ are larger than $\a$ or $\b$, respectively, for $\lm$ to be larger than $\ab$. For example, $(2)$ is clearly larger than $(1)$ but $(3,2)/(2)=\young(::\h,\h\h)$ is not larger than $(3,2)/(1)=\young(:\h\h,\h\h)$.

What are the covering relations? Let $\ab$ and $\lm$ be basic skew diagrams, then $\lm$ covers $\ab$, $\ab\co\lm$, if either 
\begin{itemize}
 \item $\lm=\ab+(1^x)/(1^y)$ with $0\leq y \leq x \leq l(\a)$ or
 \item $\lm=\ab\cup(x)/(y)$ with $0\leq y \leq x \leq \a_1$.
\end{itemize}

Note that we assumed that both $\lm$ and $\ab$ are basic. If $\lm=\ab+(1^x)/(1^y)$ with $0\leq y \leq x \leq l(\a)$ but $\lm$ is not basic, then $\lm$ does not cover $\ab$.

Note that $\lm$ has exactly one non-empty row or column more than $\ab$ if $\lm$ covers $\ab$. So the above partial order is a graded partial order with ranking function $\r(\lm)=\l_1+l(\l)$ the number of non-empty rows and columns of the basic skew diagram $\lm$.

\begin{Def}
 Let $\lm$ be a basic skew diagram, and let $n$ be minimal with $dp(\m)+1,dp(\l)\geq n$. We then say, that $\lm$ has $\d$ value $n$ and write $\d(\lm)=n$.
\end{Def}

For example $\lm=\young(::::\h\h\h,::\h\h\h\h\h,:\h\h\h,\h\h\h,\h\h,\h)$ has $\d$ value $4$. Note, that the $\d$ value of a given skew diagram is the minimal number of convex corners of either the inner or outer path. Furthermore,  the skew diagram consisting of $n$ disconnected boxes $\yng(1)$ has $\d$ value $n$, so $\d(\dn)=n$. We would like to show that $\lm$ with $\d(\lm)=n$ is larger than $\dn$, but this is false.

Take for example the skew diagram $\lm=\young(:::\h\h,:::\h,::\h,:\h,\h\h)$ with $\d(\lm)=4$. It is easy to see, that $\lm$ can not be obtained from $\d_4/\d_3$ by repeatedly applying $+$ and $\cup$ to $\d_4/\d_3$ and so $\lm$ is not larger than $\d_4/\d_3$. On the other hand $\ab=\young(::::\h,::\h\h,::\h,:\h,\h\h)$ is obtained by reordering the disconnected components of $\lm$ and we have 
\[\left( \young(:::\h,::\h,:\h,\h)+(1^4)/(1^1) \right) \cup(2)/(2)=\young(::::\h,::\h\h,:\h\h,\h\h)\cup(2)/(2)=\young(::::\h,::\h\h,::\h,:\h,\h\h)=\ab .\]

To fix this, we will now introduce an equivalence relation on the set of basic skew diagrams, so that $\overline{\lm}=\overline{\ab}$ if $\lm$ and $\ab$ are the same up to translation of the skew diagrams into which $\lm$ and $\ab$ may decompose. For example $\lm=\young(::\h,\h\h)$ and $\ab=\young(:\h\h,\h)$ both decompose into $\yng(1)\otimes\yng(2)$ and so $\overline{\lm}=\overline{\ab}$. We may now define a partial order on the set of these equivalence classes by giving the cover relations and assume transitivity. Let $\overline{\lm}$ cover $\overline{\ab}$ if there is a skew diagram  $\A\in\overline{\lm}$ which covers a skew diagram $\B\in\overline{\ab}$. All skew diagrams in the same equivalence class have the same number of non-empty rows and columns. From this follows that also the poset of equivalence classes is graded with grading function $\rho(\overline{\lm})=\l_1+l(\l)$. Notice that $\overline{\lm}>\overline{\ab}$ does not require the existence of skew diagrams $\A\in\overline{\lm},\B\in\overline{\ab}$ with $\A>\B$. Setting $\d(\overline{\lm})=\d(\lm)$ is well defined.

\begin{Le}\label{Le:lmdecays}
 Let $\lm=(1)\otimes\ab$ with $\d(\ab)=n$ and assume that $\overline{\ab}\geq \overline{\dn}$.

Then $\overline{\lm}\geq\overline{ \d_{n+1}/\d_{n}}$.
\end{Le}
\begin{proof}
We may assume, that $\lm$ has in the lower left corner the single box $(1)$ and atop to the right the skew diagram $\ab$.

Since $\overline{\ab}\geq\overline{\dn}$ there is a sequence of covering skew diagrams from $\dn$ to $\ab$. So we can choose $\circ^i\in\{+,\cup\}$ and $\A^i\in\{(1^{a+b})/(1^{b}),(a+b)/(b)\}$ such that with $\B^0=\dn, \B^m=\ab$ and $\C^i=\B^{i-1}\circ^i\A^i$ with $\overline{\C^i}=\overline{\B^i}$ we have $\overline{\B^i}\co\overline{\B^{i+1}}$.

Now, set $\widetilde{\A}=(a+b+1)/(b+1)$ if $\A=(a+b)/(b)$ and $\widetilde{\A}=\A$ otherwise.

Let  $\widetilde{\B^0}=\d_{n+1}/\d_n, \widetilde{\C^m}=\lm$ and let $\widetilde{\B^i}=(1)\otimes\B^i$ and $\widetilde{\C^i}=(1)\otimes\C^i$ such that both $\widetilde{\B^i},\widetilde{\C^i}$ contain a single disconnected box in the lower left corner. We then have $\widetilde{\C^i}=\widetilde{\B^{i-1}}\circ^i\widetilde{\A^i}$ with $\overline{\widetilde{\C^i}}=\overline{\widetilde{\B^i}}$ and so  $\overline{\widetilde{\B^i}}\co\overline{\widetilde{\B^{i+1}}}$. 
\end{proof}

\begin{Sa}\label{Sa:posetmain}
Let $\lm$ be a basic skew diagram with $\d(\lm)=n$.

Then $\overline{\lm}\geq\overline{\dn}$.
\end{Sa}
\begin{proof}
We will give a procedure to reduce $\lm$ by one rank without changing the $\d$ value. Repeatedly applying this procedure will result in a minimal skew diagram with fixed $\d$ value and we will see, that we can always reduce the rank by one without changing the $\d$ value, unless $\lm=\dn$. This shows, that $\overline{\dn}$ is the unique minimal element with $\d(\overline{\lm})=n$ and that all $\overline{\lm}$ with $\d(\overline{\lm})=n$ are larger than $\overline{\dn}$.

To give this procedure we will call a pair $(X_i,X_o)$ where $X_i$ denotes the $i$th step of the inner sequence and $X_o$ denotes the $o$th step of the outer sequence a \textit{removable pairing} if both are either $h$'s or $v$'s and if either the inner $h$ is weakly atop the outer $h$ or the inner $v$ is weakly to the left of the outer $v$. If $X_i=X_o=h$ we will call this an $h$ pairing and if both are $v$'s we will call it an $v$ pairing.

If $\ab$ is obtained from $\lm$ by removing a removable pairing then $\lm=\ab+(1^{a+b})/(1^a)$ in case of $h$ pairings and $\lm=\ab\cup(a+b)/(a)$ in case of $v$ pairings. In both cases $\lm\geq \ab$.

Because of Lemma~\ref{Le:lmdecays} we may assume, that $\lm$ does not decay into $(1)\otimes\A$ with $\A$ some arbitrary skew diagram. If $\lm$ would decay in this way, then $\d(\A)=n-1$ and if we prove that $\overline{\A}\geq\overline{ \d_{n-1}/\d_{n-2}}$ then by Lemma~\ref{Le:lmdecays} $\overline{\lm}\geq\overline{\dn}$.

We have the following possibilities for the skew diagram. We may assume that in each case none of the previous case applied.

\begin{enumerate}
 \item $\lm=\dn$ (or to be precise $\lm=\emptyset$, because we assumed $\lm\neq(1)\otimes\A$). Then there is nothing to prove.

  \item\label{2} Suppose there is a removable $h$ pairing $(h_i,h_o)$  such that both $h_i$ and $h_o$ are next to another $h$ in the inner and outer sequence, respectively, and removing it reduces the rank by one. Then the inner and outer ways are as follows
\begin{align*}
 i&: &&\ldots \ldots h h_i \ldots \\
 o&: && \ldots h h_o \ldots \ldots.
\end{align*}
Then we can remove this pair and reduce the rank by one without changing the $\d$ value. The same applies to $v$ pairings instead of $h$ pairings.

\item  Suppose there is a removable $h$ pairing $(h_i,h_o)$ such that both $h_i$ and $h_o$ are next to another $h$ in the inner and outer sequence, respectively,  but removing this pairing gives an $\ab$ which has rank more than one less than $\lm$. So the inner and outer way are as follows
\begin{align*}
 i&: &&\ldots \ldots h h_i \ldots \\
 o&: && \ldots h h_o \ldots \ldots.
\end{align*} Let $(h_i,h_o)$ be the pairing in $\lm$ such that $i-o$ is minimal of all pairings we could choose. Because $\lm$ is basic and the $h$ pairing can be removed it is $i>o$. Since $\lm$ doesn't cover $\ab$ it follows that $\ab$ can't be basic. So there has to be a  $k$ such that in $\ab$ the inner  sequence $i_1\ldots i_kX\ldots$ has in the first $k$ positions the same number of $h$'s (and $v$'s) as the outer sequence $o_1\ldots \widehat{h_o}\ldots o_{k+1}X\ldots$ has in the first $k$ positions and both continue with the same step $X\in\{h,v\}$, where $ \widehat{h_o}$ means, that $h_o$ was removed.

So  for $\lm$ we have
\begin{align*}
\lm:&& i&:  &\ldots \ldots \ldots & i_k \ldots h_i\ldots \\
&& o&: & \ldots h_o \ldots &o_k  \ldots\ldots\end{align*} 
while for $\ab$ we have
\begin{align*}
\ab:&& i&:  &\ldots \ldots \ldots & i_k X \ldots  \\
&& o&: & \ldots \widehat{h_o}\ldots &o_{k+1}  X \ldots .\end{align*} Let $k$ be minimal with this property.

Since $k$ is minimal, we have $i_k\neq o_{k+1}$. But $i_k=h$ and $o_{k+1}=v$ is not possible, because $\lm$ is basic and a skew diagram.

So we have $i_k=h$ and $o_{k+1}=v$.

If $X=h$ (so $i_{k+1}=h$) then the pairing $(i_{k+1},h_o)$ would be a removable pairing and by choice $k+1-o < i-o$ which contradicts the minimality of $i-o$. The pairing $(i_{k+1},h_o)$ would also be removable using (\ref{2}) because the minimality of $k$ assures there cannot appear non-basic configurations between the positions $o$ and $k+1$.

So we have $X=v$ (see Figure~\ref{pic1}).

\begin{figure}[h]\begin{center}
\psset{xunit=0.35cm,yunit=0.35cm,runit=0.35cm}
\begin{pspicture*}(-2,-1)(9,5)
\psline(0,2)(2,2)(2,4)
\psline(4,0)(4,4)\psdot*(4,2)
\put(-2,3){$k$}\psline{->}(-1.5,3)(1,2)
\put(-2,4){$k+1$}\psline{->}(0.5,4.25)(2,3)
\put(6,0){$k+1$}\psline{->}(6,0.25)(4,1)
\put(6,1){$k+2$}\psline{->}(6,1.25)(4,3)
\end{pspicture*}
\caption{after $X=v$ determined} \label{pic1}\end{center}
\end{figure}
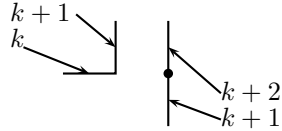

If $i_{k+2}=v$ we could remove $(i_{k+1},o_{k+2})$. This would remove one row and in this situation reduce the rank by only $1$ without changing the $\d$ value.

So we have $i_{k+2}=h$ (see Figure~\ref{pic2}).

\begin{figure}[h]\begin{center}
\psset{xunit=0.35cm,yunit=0.35cm,runit=0.35cm}
\begin{pspicture*}(-2,-1)(9,6)
\psline(0,2)(2,2)(2,4)(4,4)
\psline(4,0)(4,4)\psdot*(4,2)
\put(-2,3){$k$}\psline{->}(-1.5,3)(1,2)
\put(-2,4){$k+1$}\psline{->}(0.5,4.25)(2,3)
\put(-2,5){$k+2$}\psline{->}(0.5,5.25)(3,4)
\put(6,0){$k+1$}\psline{->}(6,0.25)(4,1)
\put(6,1){$k+2$}\psline{->}(6,1.25)(4,3)
\end{pspicture*}
\caption{$i_{k+2}=h$} \label{pic2}\end{center}
\end{figure}

Since $\lm$ is basic and the inner and outer path meet after the $k+2$nd step, it follows that $i_{k+3}\neq o_{k+3}$ and therefore we have
$i_{k+3}=v$ and $o_{k+3}=h$ (see Figure~\ref{pic3}).

\begin{figure}[h]\begin{center}
\psset{xunit=0.35cm,yunit=0.35cm,runit=0.35cm}
\begin{pspicture*}(-2,-1)(11,7)
\psline(0,2)(2,2)(2,4)(4,4)(4,6) 
\psline(4,0)(4,4)(6,4)                 
\psdot*(4,2)
\put(-2,3){$k$}\psline{->}(-1.5,3)(1,2)
\put(-2,4){$k+1$}\psline{->}(0.5,4.25)(2,3)
\put(-2,5){$k+2$}\psline{->}(0.5,5.25)(3,4)
\put(-2,6){$k+3$}\psline{->}(0.5,6.25)(4,5)
\put(8,0){$k+1$}\psline{->}(8,0.25)(4,1)
\put(8,1){$k+2$}\psline{->}(8,1.25)(4,3)
\put(8,2){$k+3$}\psline{->}(8,2.25)(5,4)
\end{pspicture*}
\caption{$i_{k+3}=v$ and $o_{k+3}=h$} \label{pic3}\end{center}
\end{figure}
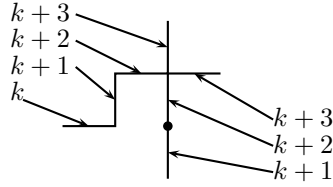

If we would have $o_{k+4}=h$, we could remove $(h_i, o_{k+4})$ which contradicts the minimality of $i-o$.

So we have $o_{k+4}=v$ (see Figure~\ref{pic4}).

\begin{figure}[h]\begin{center}
\psset{xunit=0.35cm,yunit=0.35cm,runit=0.35cm}
\begin{pspicture*}(-2,-1)(11,7)
\psline(0,2)(2,2)(2,4)(4,4)(4,6)  
\psline(4,0)(4,4)(6,4)(6,6)                
\psdot*(4,2)
\put(-2,3){$k$}\psline{->}(-1.5,3)(1,2)
\put(-2,4){$k+1$}\psline{->}(0.5,4.25)(2,3)
\put(-2,5){$k+2$}\psline{->}(0.5,5.25)(3,4)
\put(-2,6){$k+3$}\psline{->}(0.5,6.25)(4,5)
\put(8,0){$k+1$}\psline{->}(8,0.25)(4,1)
\put(8,1){$k+2$}\psline{->}(8,1.25)(4,3)
\put(8,2){$k+3$}\psline{->}(8,2.25)(5,4)
\put(8,3){$k+4$}\psline{->}(8,3.25)(6,5)
\end{pspicture*}
\caption{ $o_{k+4}=v$} \label{pic4}\end{center}
\end{figure}

But if now $i_{k+4}=h$ this would contradict $\lm\neq(1)\otimes\A$ so $i_{k+4}=v$ and we have the situation as in Figure~\ref{pic5}.

\begin{figure}[h]\begin{center}
\psset{xunit=0.35cm,yunit=0.35cm,runit=0.35cm}
\begin{pspicture*}(-2,-1)(11,8.1)
\psline(0,2)(2,2)(2,4)(4,4)(4,6)(4,8)   
\psline(4,0)(4,4)(6,4)(6,6)                  
\psdot*(4,2)\psdot*(4,6)
\put(-2,3){$k$}\psline{->}(-1.5,3)(1,2)
\put(-2,4){$k+1$}\psline{->}(0.5,4.25)(2,3)
\put(-2,5){$k+2$}\psline{->}(0.5,5.25)(3,4)
\put(-2,6){$k+3$}\psline{->}(0.5,6.25)(4,5)
\put(-2,7){$k+4$}\psline{->}(0.5,7.25)(4,7)
\put(8,0){$k+1$}\psline{->}(8,0.25)(4,1)
\put(8,1){$k+2$}\psline{->}(8,1.25)(4,3)
\put(8,2){$k+3$}\psline{->}(8,2.25)(5,4)
\put(8,3){$k+4$}\psline{->}(8,3.25)(6,5)
\end{pspicture*}
\caption{$i_{k+4}=v$} \label{pic5}\end{center}
\end{figure}

But now $(i_{k+3},o_{k+2})$ is a removable $v$ pairing and removing it changes the rank by one without altering the $\d$ value. 

The same applies to $v$ pairings instead of $h$ pairings. This means, that it is not possible to have only removable $h$ or $v$ pairing whose removal would, without altering the $\d$ value, reduce the rank by more than one.

\item Suppose now that there are only non removable pairings $(X_i,X_o)$ with $X=X_i=X_o\in\{h,v\}$ such that $X_i$ is next to another $X$ in the inner sequence and $X_o$ is next to another $X$ in the outer sequence. Suppose this is an $h$ pairing. For the skew diagram this means, that the outer $h$ is in a higher position than the inner $h$ (see Figure~\ref{pic6}).

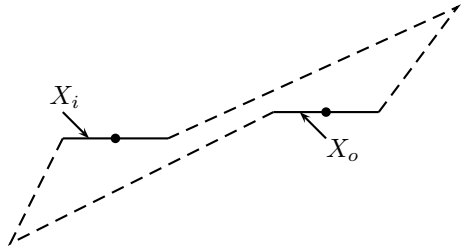
\begin{figure}[h]\begin{center}
\psset{xunit=0.35cm,yunit=0.35cm,runit=0.35cm}
\begin{pspicture*}(-.2,-1.2)(17.2,8.2)
\psline(10,4)(14,4)\psline(6,3)(2,3)
\psline[linestyle=dashed](2,3)(0,-1)(10,4)
\psline[linestyle=dashed](14,4)(17,8)(6,3)
\psdot*(12,4)\psdot(4,3)
\put(1.5,4.3){$X_i$}\psline{->}(2,4)(3,3)
\put(12,2.3){$X_o$}\psline{->}(12,3)(11,4)
\end{pspicture*}
\caption{$\lm$ has only non removable $h$ pairings} \label{pic6}\end{center}
\end{figure}

Since $\lm$ is basic, the outer sequence starts with an $h$ and because there are no removable pairings it continues with an $v$. Because $\lm$ doesn't decay into a single box and another skew diagram the inner sequence has to start with $vv$. If now the outer sequence would contain a subsequence $vv$ this would give a removable pairing, so the outer sequence does not contain a subsequence $vv$. Because $\lm$ is basic, the outer sequence ends with an $v$ and because it doesn't contain the subsequence $vv$ it ends with $hv$. Since $\lm$ doesn't  decay into a single box and another skew diagram the inner sequence has to end with $hh$. This $h$ in the inner sequence together with $X_o$ from the outer sequence form a removable pairing (see Figure~\ref{pic7}).

\begin{figure}[h]\begin{center}
\psset{xunit=0.35cm,yunit=0.35cm,runit=0.35cm}
\begin{pspicture*}(-2.2,-3.2)(19.2,10.2)
\psline(-2,1)(-2,-3)(0,-3)(0,-1)\psline(10,4)(14,4)\psline(17,8)(19,8)(19,10)(15,10)\psline(6,3)(2,3)
\psline[linestyle=dashed](2,3)(-2,1)
\psline[linestyle=dashed](0,-1)(10,4)
\psline[linestyle=dashed](14,4)(17,8)
\psline[linestyle=dashed](15,10)(6,3)
\psdot*(12,4)\psdot(4,3)\psdot*(-2,-1)\psdot*(17,10)
\put(17,4.3){$h$}\psline{->}(17,5)(16,10)
\put(12,2.3){$X_o$}\psline{->}(12,3)(11,4)
\end{pspicture*}
\caption{The removable pairing} \label{pic7}\end{center}
\end{figure}
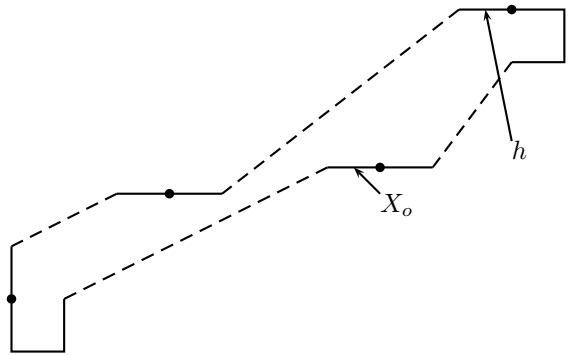

\item So we may now assume that there are no pairings $(X_i,X_o)$ with $X=X_i=X_o\in\{h,v\}$ such that $X_i$ is next to another $X$ in the inner sequence and $X_o$ is next to another $X$ in the outer sequence. By rotation symmetry we may assume that $\l$ has strictly more different parts than $\m$ (otherwise exchange the inner and outer sequence). Since there exists none of the above pairings and $dp(\l)>dp(\m)$ it follows that $\l=\d_m$ for some $m$ and, furthermore, that we have for the inner sequence either $i:\qquad \ldots h h_i\ldots$ or $i:\qquad \ldots v v_i\ldots$.

Suppose we are in the first case that we have  $i:\qquad \ldots h h_i\ldots$ (see Figure~\ref{case5}).

\begin{figure}[h]\begin{center}
\psset{xunit=0.35cm,yunit=0.35cm,runit=0.35cm}
\begin{pspicture*}(-2.1,-2.1)(12.1,14.1)
\psline(0,0)(0,2)(2,2)(2,4)(4,4)(4,6)(6,6)(6,8)(8,8)(8,10)(10,10)(10,12)
\psline[linestyle=dashed](-2,-2)(0,0)\psline[linestyle=dashed](10,12)(12,14)
\psline(0,12)(4,12)\psdot*(2,12)\psline[linestyle=dashed](-2,10)(0,12)\psline[linestyle=dashed](4,12)(6,14)
\put(2.5,12.25){$h_i$}
\end{pspicture*}
\caption{$h_i$ next to another $h$} \label{case5}\end{center}
\end{figure}
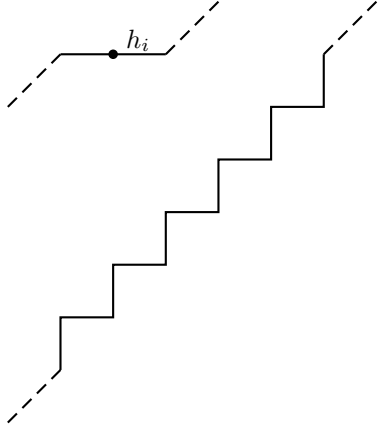

Then we can remove the column containing $h_i$ and by doing so reduce the rank by one without changing the $\d$ value (see Figure~\ref{case52}).

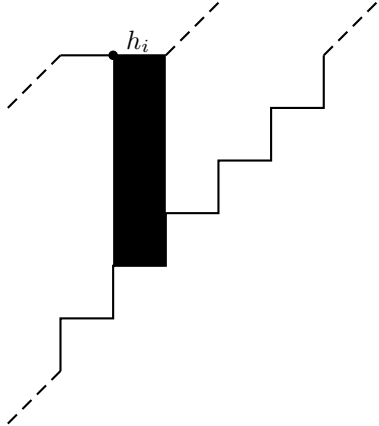
\begin{figure}[h]\begin{center}
\psset{xunit=0.35cm,yunit=0.35cm,runit=0.35cm}
\begin{pspicture*}(-2.1,-2.1)(12.1,14.1)
\psline(0,0)(0,2)(2,2)(2,4)(4,4)(4,6)(6,6)(6,8)(8,8)(8,10)(10,10)(10,12)
\psline[linestyle=dashed](-2,-2)(0,0)\psline[linestyle=dashed](10,12)(12,14)
\psline(0,12)(4,12)\psdot*(2,12)\psline[linestyle=dashed](-2,10)(0,12)\psline[linestyle=dashed](4,12)(6,14)
\put(2.5,12.25){$h_i$}
\psframe*(2,12)(4,4)
\end{pspicture*}
\caption{The removable pairing} \label{case52}\end{center}
\end{figure}
\end{enumerate}

This finishes the proof.
\end{proof}

\section{Application to skew characters: \\Lower bounds for the number of components, constituents and pairs of components which differ by one box}\label{sec:fcmain:bounds}

In this section we are interested in skew characters and so do not strictly distinguish between the skew diagrams and equivalence classes of skew diagrams up to translation.

\begin{Def}
We say that a skew diagram $\A$ or skew character $[\A]=[\lm]=\sum_\n c(\l;\m,\n)$ is of cc-type $(a,b)$ if $[\A]$ has $a=\sum_{c(\l;\m,\n)\neq 0}1$ components and $b=\sum_\n c(\l;\m,\n)$ constituents. We then also write $cc(\A)=(a,b)$ or $cc([\A])=(a,b)$. Note that always $a \leq b$ so there is no way of confusing the order. Furthermore we say that $\A$ with $cc(\A)=(a,b)$ has cc-type at least $(c,d)$ if $a\geq c$ and $b\geq d$.
\end{Def}

For example, the skew character corresponding to $(2,1)\otimes(2,1)=\young(::\h\h,::\h,\h\h,\h)$ is
\[[(2,1)]\otimes[(2,1)]=[4,2]+[4,1^2]+[3^2]+2[3,2,1]+[3,1^3]+[2^3]+[2^2,1^2]\] 
and so $cc((2,1)\otimes(2,1))=cc([(4,3,2,1)/(2^2)])=(7,8)$.

For the following proofs we use the following lemma which we proved in~\cite{GSLR} and is a generalization of a lemma in~\cite{Gut}.
\begin{Le}[Lemma 3.1, \cite{GSLR}]\label{Le:add}
Let $\l,\m,\n,\l',\m',\n' $ be partitions with $c(\l';\m',\n')\neq0$.

Then
\[c(\l;\m,\n) \leq c(\l+\l';\m+\m',\n+\n')\]
and by conjugation
\[c(\l;\m,\n) \leq c(\l\cup\l';\m\cup\m',\n\cup\n').\]
\end{Le}

\begin{Bem}\label{Bem:addcc-typelr}
Note that $\l^1+\n\neq\l^2+\n$ for $\l^1\neq\l^2$ so this lemma tells us that adding a skew diagram $\B$ to a skew diagram $\A$ weakly increases the number of components and constituents of $[\A+\B]$ compared to $[\A]$ (or $[\B]$). By conjugation the same applies to the row wise addition  of two skew diagrams $\A\cup\B$. This allows us to consider small examples of $[\A]$ to give a lower bound on the number of components and constituents of larger $[\A']$ if $\A'$ can be obtained from $\A$ by successively adding, column or row wise, $\B^i$ for some skew diagrams $\B^i$.
\end{Bem}

We will now introduce a partial order on the set of skew characters by giving the cover relations. Let $\chi$ and $\psi$ be skew characters, then we say that $\chi$ covers $\psi$ if there exists skew diagrams $\lm$ and $\ab$ with $\chi=[\lm]$ and $\psi=[\ab]$ such that $\lm$ covers $\ab$. Since the number of non empty rows and columns of a skew diagram is fixed for a given skew character this partial order of skew characters is also graded with ranking function $\rho([\lm])=\l_1+l(\l)$ for basic skew diagrams $\lm$. Note that this partial order is compatible with the partial order on the set of equivalence classes of skew diagrams of Section~\ref{sec:poset}.

\begin{Sa}\label{Sa:dpcc}
 Let $\lm$ be a basic skew diagram with $\d(\lm)=n$. Then $cc(\lm)$ is at least $(p_n, f_n)$ where $p_n$ is the number of partitions of $n$ and $f_n$ the number of standard Young tableaux with $n$ boxes.
\end{Sa}
\begin{proof}
Let $\d_n=(n,n-1,n-2,\ldots,2,1)$ then as an easy consequence of the LR rule we have
\[ [\d_n/\d_{n-1}]=[\underbrace{(1)\otimes(1)\otimes\cdots\otimes(1)}_{n\textnormal{-times}}]=[1]^n=\sum_{\l\vdash n} f^\l[\l] \]
where $f^\l$ is the number of standard Young tableaux of shape $\l$. So we have $cc(\d_n/\d_{n-1})=(p_n,f_n)$. Since $\d(\lm)=n$  $\lm$ is
larger than $\d_n/\d_{n-1}$ by Theorem~\ref{Sa:posetmain} and so $cc(\lm)$ is at least $(p_n,f_n)$.
\end{proof}

We will use the following notation  in the remaining part of this chapter.
\begin{Def}
 
We let $\bp_n$ denote the number of partitions of $n$ with two different kinds of $1$'s and $2$'s. For the partitions of $2$ with two different kinds of  $1$'s and $2$'s see Example~\ref{Bsp:dp}.

Let $g_n$ denote  the number of unordered pairs $(\n^1,\n^2)$ of partitions of $n$ with $\abs{\n^1\cap\n^2}=n-1$. So $g_n$ counts the pairs of partitions of $n$ which differ only by one box.

\end{Def}

\begin{Le}\label{Le:bp_gn}
Then $\bp_n=g_{n+2}$ for all $n$.
\end{Le}
\begin{proof}
We give a bijection of partitions of $n$ with two different kinds of $1$'s and $2$'s to pairs $(\n^1,\n^2)$ of partitions of $n+2$ which differ only by one box. We may assume that $\n^1$ is lexicographically larger than $\n^2$.

Suppose the two kinds of $1$'s are the usual $1$ and the other be $1'$ and the two kinds of $2$'s are $2$ and $2'$.
Let $\bl$ be such a partition of $n$ and let $\l$ denote the partition formed by the usual parts of $\bl$.  Furthermore, let $n_1$ denote  the number of $1'$ in $\bl$ and $n_2$ denote the number of $2'$ in $\bl$. So $\bl=\l\cup({2'}^{n_2},{1'}^{n_1})$.

For a partition $\bl$ now define the bijection by setting 
\[\n^1=\l\cup(n_1+n_2+2,n_2), \qquad \n^2=\l\cup(n_1+n_2+1,n_2+1). \]

Now obviously $\n^1$ is lexicographically larger than $\n^2$ and both partitions differ only by one box. Furthermore, different $\bl$ correspond to different triples $(\l,n_1,n_2)$ and so give different pairs $(\n^1,\n^2)$.

Finally the inverse map is obtained as follows. If $\n^1$ and $\n^2$ differ by only one box (and $\n^1$ is lexicographically larger than $\n^2$) , then $\n^2$ is obtained from $\n^1$ by removing a box in one row and placing it in a lower row. Let all the other rows form $\l$ then the two rows which are different are of the form $(a+1)$ and $(b)$ in $\n^1$ and $(a)$ and $(b+1)$ in $\n^2$ for $a\geq b\geq 0$. Now $a+1>b+1$ since otherwise $\n^1=\n^2$. So to exclude this case we may instead assume that the rows are $(c+2)$ and $(b)$ in $\n^1$ and $(c+1)$ and $(b+1)$ in $\n^2$ for $c\geq b\geq 0$.
Setting $n_1=c-b$ and $n_2=b$ gives the inverse map.
\end{proof}

\begin{Bsp}\label{Bsp:dp}
 We have $\bp_2=5$ and there is the following correspondence given by the above bijection.

 \begin{align*}
 \bl   && \l     & &          \n^1     && \n^2                 \\
    2   && 2       & &  \yng(2,2)    && \yng (2,1,1)     \\
   2'   && 0       & &  \yng(3,1)    && \yng (2,2)        \\
   1^2  && 1^2  & &  \yng(2,1,1) && \yng (1,1,1,1) \\
1,1'  &&  1   & &  \yng(3,1)    && \yng (2,1,1)      \\
{1'}^2 && 0    & &  \yng(4)       && \yng (3,1)
 \end{align*}
\end{Bsp}

\begin{Bem}
 Lemma~\ref{Le:bp_gn} is useful because one sees directly that the generating function for $\bp_n$ is given by 
\[ \sum_{i\geq0} \bp_i x^i= \frac{1}{(1-x)(1-x^2)}\prod_{i\geq 1} \frac{1}{1-x^i}.\]
\end{Bem}

In the following theorem the condition $\d(\lm)\geq2$ only makes sure that $\lm$ is neither a partition nor a rotated partition but constrains $\lm$ not in any other way. The case that $\lm$ is a partition $\a$ or rotated partition $\a^\circ$ is uninteresting for the theorem because then $[\lm]=[\a]$ is irreducible.

\begin{Sa}\label{Sa:2chdiffer}
Let $\lm$ be a basic skew diagram with $\d(\lm)=n\geq2$.

Then $[\lm]=\sum_\n c(\l;\m,\n)[\n]$ contains at least $g_n$ pairs of characters $([\n^1],[\n^2])$ whose corresponding diagrams differ only by one box, i.e.\ there are $\n^1,\n^2$ with $\abs{\n^1\cap\n^2}=\abs{\n^1}-1=\abs{\n^2}-1$ and $c(\l;\m,\n^1),c(\l;\m,\n^2)\neq0$ (with $g_n$ as in Lemma~\ref{Le:bp_gn}).

Furthermore, if $\l=(\l_1,\ldots,\l_l),\m=(\m_1,\ldots,\m_m)$ with $\l_l,\m_m\geq1$ set $\A=(\l_1-2,\l_l-1)/(\m_1-1)$ and $\B=(\l_2,\l_3,\ldots\l_{l-1})/(\m_2,\m_3,\ldots,\m_m)$ with $[\A]$ having $a$ components and $[\B]$ having $b$ components. Then there are at least $\max(a,b)$ of those pairs $\n^1,\n^2$.
\end{Sa}
\begin{proof}
We first show there are at least $\max(a,b)$ pairs $\n^1,\n^2$.

We can deduce this part of the theorem from the fact that $[(2,1)/(1)]=[2]+[1^2]$ contains two characters whose corresponding diagrams differ only by one box.

We explicitly show how to obtain $\lm$ from $(2,1)/(1)$.

The skew diagram $(\l_1,\l_l)/(\m_1)$ is larger than $(2,1)/(1)$ 
\begin{align*}
(\l_1,\l_l)&=(2,1)+(\l_1-2,\l_l-1), &(\m_1)&= (1)+(\m_1-1)
\end{align*}
and $\A=(\l_1-2,\l_l-1)/(\m_1-1)$ is a skew diagram. Let $\a$ be a partition such that $[\a]$ appears in $[\A]$, so $c((\l_1-2,\l_l-1);(\m_1-1),\a)\neq 0$.

Then by Lemma~\ref{Le:add} $[\a+(1^2)]$ and $[\a+(2)]$ both appear in $[(\l_1,\l_l)/(\m_1)]$ and, furthermore, $\a+(1^2)\cap\a+(2)=\a+(1)$ so $\a+(1^2)$ and $\a+(2)$ differ by only one box.

Now $\lm$ is larger than $(\l_1,\l_l)/(\m_1)$ 
\begin{align*}
\l&=(\l_1,\l_l)\cup(\l_2,\l_3,\ldots,\l_{l-1}),& \m&= (\m_1)\cup(\m_2,\m_3,\ldots,\m_m)
\end{align*}
and $\B=(\l_2,\l_3,\ldots\l_{l-1})/(\m_2,\m_3,\ldots,\m_m)$ is a skew diagram. Let $\b$ be a partition such that $[\b]$ appears in $[\B]$.

Then by Lemma~\ref{Le:add} $[(\a+(1^2))\cup\b]$ and $[(\a+(2))\cup\b]$ both appear in $[\lm]$ and
\[\left((\a+(1^2))\cup\b\right)\cap\left((\a+(2))\cup\b\right)=(\a+(1))\cup\b\]
so $\n^1=(\a+(1^2))\cup\b$ and $\n^2=(\a+(2))\cup\b$ differ only by one box.

Furthermore, notice  that a different choice for $\a$ or $\b$ yields a different pair $\n^1,\n^2$. This proves that there are at least $\max(a,b)$ pairs $\n^1,\n^2$.

Now we will prove that there are also at least $g_n$ pairs $\n^1,\n^2$.

As mentioned above, as an easy consequence of the LR rule we have
\[ [\d_n/\d_{n-1}]=[\underbrace{(1)\otimes(1)\otimes\cdots\otimes(1)}_{n\textnormal{-times}}]=[1]^n=\sum_{\l\vdash n} f^\l[\l] \]
where $f^\l$ is the number of standard Young tableaux of shape $\l$, in particular, all irreducible characters of $S_n$ appear in $[\d_n/\d_{n-1}]$. So by definition of $g_n$ $[\d_n/\d_{n-1}]$ contains $g_n$ characters $[\a],[\b]$ whose corresponding diagrams differ only by one box.

By Theorem~\ref{Sa:posetmain} $\lm$ is larger than $\d_n/\d_{n-1}$, so there exist skew diagrams $\B^i$ such that $\lm$ is obtained from $\dn$ by using the operations $+,\cup$ together with the $\B^i$. Let $\circ^i$ be either $+$ or $\cup$ then 
\[ \lm=((\dn\circ^1\B^1)\circ^2\B^2)\cdots\circ^j\B^j.\]
 Choose $[\a^i]$ contained in $[\B^i]$ and $[\bn^1],[\bn^2]$ contained in $[\dn]$ with $\abs{\bn^1\cap\bn^2}=n-1$. Set
\begin{align*}
 \n^1&=((\bn^1\circ^1\a^1)\circ^2\a^2)\cdots\circ^j\a^j, &\n^2&=((\bn^2\circ^1\a^1)\circ^2\a^2)\cdots\circ^j\a^j
\end{align*}
then by Lemma~\ref{Le:add} both $[\n^1],[\n^2]$ appear in $[\lm]$ and, furthermore, $\abs{\n^1\cap\n^2}=\abs{\n^1}-1$. Finally a different choice of $\bn^1,\bn^2$ gives different $\n^1,\n^2$ (for fixed $(\a^i,\circ^i)$) and there are by definition $g_n$ choices for $\bn^1,\bn^2$.
\end{proof}

\begin{Bem}\label{Bem:firstgpf}
 In the On-Line Encyclopedia of Integer Sequences~\cite{OEIS} $g_{n}=\bp_{n-2}$ has the id: A000097, $p_n$ has the id: A000041 and $f_n$ has the id: A000085. Their first terms are
\[
\begin{array}{l|ccc|ccc|ccc|ccc|c}
n:  &1 &2 &3 &4 &5 &6 &7 &8 &9 &10 &11 &12 &13 \\
\hline
g_n:  &0 & 1 &2 &5 &9 &17 &28 & 47&73 & 114 &170 &253 &365 \\
p_n: &1 & 2 &3 &5 & 7 &11 &15 &22  &30  &42 &56  &77  &101  \\
f_n:   &1 &2 &4 &10 &26 &76 &232 &764 &2620 &9496 &35696 &140152 &568504
\end{array}\]
\end{Bem}

\begin{Le}\label{Le:prodccpairs}
Let $\a,\b$ be partitions with $dp(\a)\geq dp(\b)=n$. Then $[\a]\otimes[\b]$ has cc-type at least $(p_{n+1},f_{n+1})$ and contains $g_{n+1}$ pairs of components $([\n^1],[\n^2])$ such that their corresponding partitions differ only by one box.
\end{Le}
\begin{proof}
This follows directly from the previous theorems by setting $\lm=\a \otimes \b^\circ$ because then $dp(\l)=dp(\a)+1,dp(\m)=dp(\b)$.
\end{proof}

\begin{Le}
Let $\lm$ be a skew diagram with $\abs{\lm}=n$.

Then $[\lm]$ contains at most
\begin{itemize}
 \item $g_n$ pairs $([\n^1],[\n^2])$ such that $\abs{\n^1\cap\n^2}=n-1$,
 \item $p_n$ components,
 \item $\min (f_n,p_nf^\m,p_nf^\bl)$ constituents (with $\bl=(\l_1-\l_l,\l_1-\l_{l-1},\ldots,\l_1-\l_3,\l_1-\l_2,0)$).
\end{itemize}
\end{Le}
\begin{proof}
 The first two statements are trivial, because there are not more irreducible characters of $S_n$.

For the third statement notice, that $\lm$ is smaller than $\dn$ which gives by Lemma~\ref{Le:add} $c(\l;\m,\n)\leq c(\d_n;\d_{n-1},\n)=f^\n$.
Since the LR coefficient is symmetric in $\m$ and $\n$ we also have $c(\l;\m,\n)\leq f^\m$ and by rotation symmetry $c(\l;\m,\n)\leq f^\bl$.

So  for the number of constituents of $[\lm]$ 
\begin{align*}
 \sum_\n c(\l;\m,\n) &\leq \sum_\n f^\n=f_n, \\
 \sum_\n c(\l;\m,\n) &= \sum_{\n\vdash n} c(\l;\m,\n) \leq \sum_{\n\vdash n} f^\m =p_nf^\m,\\
 \sum_\n c(\l;\m,\n) &= \sum_{\n\vdash n} c(\l;\m,\n) \leq \sum_{\n\vdash n} f^\bl =p_nf^\bl.
\end{align*}

Notice that all three bounds are reached for $\lm=\dn$.
\end{proof}

{\bfseries Acknowledgement:}  I want to thank Martin Rubey and Christine Bessenrodt for helpful discussions. John Stembridge's "SF-package for maple" \cite{stemmaple} was very helpful for computing examples.The superseeker service of the OEIS \cite{OEIS} was helpful in finding the sequence $\bp(n)$ of Lemma~\ref{Le:bp_gn} after the sequence $g_n$ was calculated.

\end{document}